\documentclass[12pt]{amsart}
\usepackage{amsmath, amssymb, amsfonts, enumerate}
\usepackage[all]{xy}
\usepackage{amscd}

\newtheorem{theorem}{Theorem}[section]
\newtheorem{lemma}[theorem]{Lemma}

 \theoremstyle{definition}
 
 \newtheorem{remark}[theorem]{Remark}

\numberwithin{equation}{section}

\newcommand {\N}{\mathbb{N}} 
\newcommand {\R}{\mathbb{R}} 

\DeclareMathOperator{\Id}{Id}

\begin{document}
 \title[The Surjectivity of the combinatorial Laplacian]{The surjectivity of the combinatorial Laplacian on infinite graphs}
 \author[T. Ceccherini-Silberstein]{Tullio Ceccherini-Silberstein}
\address{Dipartimento di Ingegneria, Universit\`a del Sannio, C.so
Garibaldi 107, 82100 Benevento, Italy}
\email{tceccher@mat.uniroma1.it}
\author[M. Coornaert]{Michel Coornaert}
\address{Institut de Recherche Math\'ematique Avanc\'ee,
UMR 7501,                                             Universit\'e  de Strasbourg et CNRS,
                                                 7 rue Ren\'e-Descartes,
                                               67000 Strasbourg, France}
\email{coornaert@math.unistra.fr}
\author[J. Dodziuk]{J\'ozef Dodziuk}
\address{Ph.D.\ Program in Mathematics, Graduate Center (CUNY), New York, NY 10016, USA}
\email{jdodziuk@gc.cuny.edu}
\date{\today}
\begin{abstract}
Given a connected locally finite simplicial graph $ G$ with vertex set $V$, the combinatorial Laplacian 
$\Delta_G \colon \R^V \to \R^V$ is defined on the space of all real-valued functions on  $V$.
  We prove that $\Delta_G$ is surjective if $G$ is infinite. 
 \end{abstract}
\maketitle

\section{Introduction}

Let $G $ be a connected locally finite  graph  with vertex set $V$.
To simplify the exposition, we shall always assume that $G$ is simplicial, that is, without loops and multiple edges.
Two vertices $v,w \in V$ are called \emph{adjacent}, and one then writes $v \sim w$,  if $\{v,w\}$ is an edge of $G$.

 The \emph{combinatorial Laplacian} on $G$ is the linear map $\Delta_G \colon \R^V \to \R^V$,
 where $\R^V$ is the vector space consisting of all real-valued functions on $V$, defined by
$$
\Delta_G(f)(v) = f(v) - \frac{1}{\deg(v)} \sum_{v \sim w} f(w)
$$
for all $f \in \R^V$ and $v \in V$. 
Here $\deg(v)$ denotes the \emph{degree} of the vertex $v$, i.e., the number of vertices in $G$ which are adjacent to $v$.
\par
Note that $\Delta_G$ is never injective since all constant functions are in its kernel.
 As a consequence, when the graph $G$ is finite,  $\Delta_G$ is not surjective since, in this case, $\R^V$ is 
finite-dimensional.
In the present paper, we shall establish the following result.

\begin{theorem}
\label{t:lapl-surj}
Let $G $ be an infinite, connected, locally finite simplicial  graph with vertex set $V$.
Then the combinatorial Laplacian
$\Delta_G \colon \R^V \to \R^V$ is surjective.
\end{theorem}

In the particular case when $G$ is the Cayley graph of an infinite  finitely generated group, this result had been previously obtained in \cite{laplace} by using two steps.
The first one consisted in showing that the image of the Laplacian is closed in the prodiscrete topology (see below). The second one distinguished two cases according to whether the group was amenable or not.
The proof we present here for general graphs is simpler in the sense that we also first establish the closed image property of the Laplacian (Section \ref{sec:closed}) but do not need to introduce amenability considerations. 
Instead, we apply the maximum principle to finitely-supported functions on vertices in order to prove that the image of the Laplacian is also dense in the prodiscrete topology (Section \ref{sec:surjectivity}). The image, being both closed and dense, 
must be  equal to the whole space $\R^V$.

\section{The closed image property}
\label{sec:closed}

Let $G $ be a connected locally finite simplicial graph with vertex set $V$.

The \emph{prodiscrete topology} on $\R^V$ is the product topology obtained by taking the discrete topology on each factor $\R$ of $\R^V$.
This topology is metrizable.
Indeed, if $(\Omega_n)_{n \in \N}$ is a non-decreasing sequence of finite subsets of $V$
whose union is $V$, then the metric $\delta$ on $\R^V$ defined by
$$
\delta(f,g) = \sum_{n \in \N} \frac{1}{2^{n + 1}} \delta_n(f,g) \quad \text{for all } f,g \in \R^V,
$$
where $\delta_n(f,g) = 0$ if $f$ and $g$ coincide on $\Omega_n$ and $\delta_n(f,g) = 1$ otherwise, induces the prodiscrete topology on $\R^V$.
Note that   a base of neighborhoods of $f \in \R^V$ in the prodiscrete topology is provided by the sets
$$
W_n(f) = \{g \in \R^V : f\vert_{\Omega_n} = g\vert_{\Omega_n} \}.
$$

The goal of this section is to prove that the image of $\Delta_G$ is closed in $\R^V$ in the prodiscrete topology (Lemma \ref{l:closed-image}). The proof is analogous to the proof  of the closed image property for linear cellular automata over groups whose alphabets are   finite-dimensional vector spaces 
(see \cite{gromov-esav}, \cite{ca-and-groups-springer}, \cite{periodic}). 
  \par
For completeness,  let us first recall some elementary facts about projective sequences and the Mittag-Leffler condition (cf. \cite{periodic}).
\par
A \emph{projective sequence} of sets  consists of a sequence $(X_n)_{n \in \N}$ of sets together with maps $u_{nm} \colon X_m \to X_n$ defined for all $n \leq m$
  satisfying the following conditions:

\begin{enumerate}[(PS-1)]
\item
$u_{n n}$ is the identity map on $X_n$ for all $n \in \N$;
\item
$u_{n k} = u_{n m} \circ u_{m k}$ for all $n,m,k \in \N$ such that $n \leq m \leq k $.
\end{enumerate}
Such a projective sequence will be denoted by $(X_n,u_{n m})$ or simply $(X_n)$.
\par
The \emph{projective limit} $\varprojlim X_n$ of the projective sequence $(X_n,u_{n m})$ is the subset of $\prod_{n \in \N} X_n$ consisting of all the sequences $(x_n)_{n \in \N}$ which satisfy 
$x_n = u_{n m}(x_m)$ for all $n,m \in \N$ with $n \leq m$.
\par
Observe that if $(X_n,u_{n m})$ is a projective sequence of sets then it follows from (PS-2) that, for $n \in \N$ fixed, the sequence $(u_{n m}(X_m))_{m \geq n}$ is a non-increasing sequence of subsets of $X_n$.
We say that the projective sequence $(X_n,u_{n m})$ satisfies the \emph{Mittag-Leffler condition} if,
for each $n \in \N$,  the sequence
$(u_{n m}(X_m))_{m \geq n}$ stabilizes, that is,  there exists an integer  $m_0 = m_0(n) \geq n$ such that $u_{n m}(X_m) = u_{n m_0}(X_{m_0})$ for all $m \geq m_0$.

\begin{lemma}[Mittag-Leffler]
\label{l;ML}
If  $(X_n,u_{n m})$ is a projective system of nonempty sets
which satisfies the Mittag-Leffler condition,
then its projective limit $\varprojlim X_n$ is not empty.
\end{lemma}

\begin{proof}
Let $(X_n,u_{n m})$ be an arbitrary projective sequence of sets.
 The set $X_n' = \bigcap_{m \geq n} u_{nm}(X_m)$ is called the set of \emph{universal elements} in $X_n$ (cf. \cite{grothendieck-ega-3}).
It is clear that the map $u_{n m}$ induces by restriction a map $u_{n m}' \colon X_m' \to X_n'$ for all $n \leq m$ and  that $(X_n',u_{n m}')$ is a projective sequence 
having the same projective limit as the projective sequence $(X_n,u_{n m})$.
\par
Suppose now that all the sets $X_n$ are nonempty and that the projective sequence $(X_n,u_{n m})$ satisfies the Mittag-Leffler condition.
Then, for each $n \in \N$, there is an integer $m_0 = m_0(n) \geq n$ such that $u_{n m}(X_m) = u_{n m_0}(X_{m_0})$ for all $m \geq m_0$.
It follows that $X_n' = u_{n m_0}(X_{m_0})$ so that, in particular, the set $X_n'$ is not empty.
We claim that the map $u_{n, n + 1}' \colon X_{n + 1}' \to X_n'$ is surjective for every $n \in \N$. 
Indeed, let us fix $n \in \N$ and suppose that $x_n' \in X_n'$.
 By the Mittag-Leffler condition, we can find an integer $p \geq n + 1$ such that
$u_{n k}(X_k) = u_{n p}(X_p)$ and $u_{n + 1,  k}(X_k) = u_{n + 1,p}(X_p)$ for all $k \geq p$. 
It follows that $X_n' = u_{n p}(X_p)$ and $X_{n + 1}' = u_{n + 1, p}(X_p)$. 
Consequently, we can find $x_p \in X_p$ such that $x_n' = u_{n p}(x_p)$.
Setting $x_{n + 1}' = u_{n + 1,p}(x_p)$, we have  $x_{n + 1}' \in X_{n + 1}'$ and 
$$
u_{n, n + 1}'(x_{n +1}') = u_{n, n + 1}(x_{n +1}') = u_{n, n + 1} \circ u_{n+1,p} (x_p) = u_{n p}(x_p) = x_n'.
$$ 
This proves  that $u_{n, n + 1}'$ is onto.
Now, as the sets $X_n'$ are nonempty,   we can construct by induction a sequence $(x_n')_{n \in \N}$ such that $x_n' = u_{n, n + 1}'(x_{n +1}') $ for all $n \in \N$.
This sequence is in the projective limit $\varprojlim X_n' = \varprojlim X_n$. This shows that $\varprojlim X_n$ is not empty.    
\end{proof}

\begin{lemma}
\label{l:closed-image}
Let $G$ be a connected locally finite simplicial graph with vertex set $V$.
Then the image of the combinatorial Laplacian
$\Delta_G \colon \R^V \to \R^V$ is closed in $\R^V$ in the prodiscrete topology. 
\end{lemma}

\begin{proof}
Let us fix a vertex $v_0 \in V$.
For each  $n \in \N$, let $B_n = \{v \in V : d_G(v_0,v) \leq n \}$ denote the closed ball of radius $n$ centered at $v_0$ with respect to the graph metric $d_G$ on $V$.
Observe that $\Delta_G$ induces by restriction a linear map
\begin{equation}
\Delta_G^{(n)} \colon \R^{B_{n + 1}} \to \R^{B_n} \label{restr-lapl}
\end{equation}
for every $n \in \N$.
\par
Suppose  that $g \in \R^V$ is in the closure of $\Delta_G(\R^V)$ in the prodiscrete topology.
Then, for each $n \in \N$, there exists $f_n \in \R^{V}$ such that
$g$ and $\Delta_G(f_n)$ coincide on $B_n$.
Consider, for each $n \in \N$, the affine
subspace $X_n \subset \R^{B_{n + 1}}$ defined by 
$$
X_n = (\Delta_G^{(n)})^{-1}(g\vert_{B_n}).
$$
Observe that $X_n \not= \varnothing$ since $f_n \vert_{B_{n + 1}} \in X_n$.
Now, for all $n \leq m$, the restriction map $\R^{B_{m + 1}} \to \R^{B_{n + 1}}$ induces an affine map
$u_{nm} \colon X_m \to X_n$. Conditions (PS-1) and (PS-2) are trivially satisfied so that
$(X_n,u_{nm})$ is a projective sequence. We claim that this projective sequence
satisfies the Mittag-Leffler condition. 
Indeed, 
 for $n$ fixed, 
 as the sequence  $u_{nm}(X_m)$, where  $m = n,n+1,\dots$, is a
non-increasing sequence of  affine subspaces of the finite-dimensional vector space $\R^{B_{n + 1}}$, 
it must
stabilize.
 It follows from Lemma \ref{l;ML} that the projective limit $\varprojlim X_n$ is nonempty.
Choose an element $(x_n)_{n \in \N} \in \varprojlim X_n$. We have that $x_n \in \R^{B_{n + 1}}$. Moreover,  
$x_{n + 1}$ coincides with $x_n$ on $B_{n + 1}$   for all $n \in \N$. 
As $V = \cup_{n \in \N} B_{n + 1}$, there exists a (unique)   
$f \in \R^V$ such that $f\vert_{B_{n + 1}} = x_n$ for all $n$. 
We have $(\Delta_G(f))\vert_{B_n}= \Delta_G^{(n)}(x_n) = g\vert_{B_n}$ for all $n$ since $x_n \in X_n$. 
Since $V = \cup_{n \in \N} B_n$, it follows that $\Delta_G(f) = g$.
This shows that $\Delta(\R^V)$ is closed in $\R^V$ in the prodiscrete topology.
 \end{proof}

\section{Surjectivity}
\label{sec:surjectivity}

We now prove Theorem \ref{t:lapl-surj}.
\par
Suppose that $G$ is an infinite, locally finite, connected simplicial graph.
We keep the notation introduced in the proof of Lemma \ref{l:closed-image}.
Let
$F_n$ denote the vector subspace of $\R^V$ consisting of all functions $f \in \R^V$ whose support is contained in $B_n$.
Consider the linear map
 $u_n \colon F_n \to F_n$ defined by
 $$
 u_n(f)(v) =
 \begin{cases}
 \Delta_G(f)(v)& \text{ if  } v \in B_n \\
 0 & \text{otherwise},
 \end{cases}
 $$
 for all $f \in F_n$ and $v \in V$.
  By the maximum principle, $u_n  $ is injective.
  Indeed, if $f \in F_n$ satisfies $u_n(f) = 0$, then we have
  $$
  \vert f(v) \vert = \left\vert \frac{1}{\deg(v)} \sum_{v \sim w} f(w) \right\vert \leq \frac{1}{\deg(v)} \sum_{v \sim w} \vert f(w) \vert,
  $$
  for all $v \in B_n$. This implies that if $v \in B_n$ satisfies $\vert f(v) \vert = M$, where $M = \max \vert f \vert$,
  then $\vert f(w) \vert = M$ for all $w \in V$ with $v \sim w$. Therefore $\vert f \vert$ is constant on $B_{n + 1}$.
  As $G$ is infinite, there are points in $B_{n + 1}$ that are not in $B_n$.
Consequently,  $f$ is identically zero.
Now the injectivity of $u_n$ implies its surjectivity since $F_n$ is finite-dimensional.
   It follows that for all $g \in \R^V$ and $n \in \N$,
  we can find $f \in F_n$ such that $\Delta_G(f)$ coincides with $g$ on $B_n$.
  This shows that $\Delta_G(\R^V)$ is dense in $\R^V$ in the prodiscrete topology.
  As $\Delta_G(\R^V)$ is also closed by Lemma \ref{l:closed-image},
  we conclude that $\Delta_G(\R^V) = \R^V$. 
This completes the proof of Theorem \ref{t:lapl-surj}.

\begin{remark}
More generally, the same proof yields, \emph{mutatis mutandis}, the surjectivity of 
$L = \Delta_G + \lambda \Id \colon \R^V \to \R^V$ for every infinite, locally finite, connected simplicial graph $G$ and any  function $\lambda \colon V \to [0,+ \infty)$ defined on the vertex set $V$ of $G$ 
(here $\Id$ is the identity map on $\R^V$).
Indeed, $L$ is linear and, for $f \in \R^V$ and $v \in V$, from  $L(f)(v) = 0$ we deduce that
  \begin{align*}
  \vert f(v) \vert &= \left\vert \frac{1}{(1 + \lambda(v))\deg(v)} \sum_{v \sim w} f(w) \right\vert \\
  &\leq \frac{1}{(1 + \lambda(v))\deg(v)} \sum_{v \sim w} \vert f(w) \vert 
  \leq \frac{1}{\deg(v)} \sum_{v \sim w} \vert f(w) \vert,
  \end{align*}
  so that the maximum principle is also satisfied by $L$.
 \end{remark}


\end{document}